\documentclass[letterpaper, 10 pt, conference]{ieeeconf}  

\IEEEoverridecommandlockouts                              
\usepackage{graphics} 
\usepackage{epsfig} 
\usepackage{amsmath} 
\usepackage{amssymb}  
\usepackage[dvipsnames]{xcolor}

\usepackage{stmaryrd} 
\usepackage{bbm} 
\usepackage{algorithmic}
\usepackage[linesnumbered,ruled,vlined]{algorithm2e}
\usepackage{url}
\usepackage{caption}
\usepackage{subcaption}
\usepackage{booktabs}
\usepackage{dsfont}
\usepackage{mathtools}

\newcommand\cX{\mathcal X}
\newcommand\cY{\mathcal Y}

\newcommand{\bhatw}{\bar{\hat{w}}}

\newcommand\bX{\boldsymbol{X}}
\newcommand\bY{\boldsymbol{Y}}

\newcommand\bZ{\boldsymbol{Z}}

\newcommand\bw{\boldsymbol{w}}

\newcommand\EE{\mathbb E}

\newcommand\RR{\mathbb R}

\newtheorem{theorem}{Theorem}[section]

\newtheorem{definition}[theorem]{Definition}

\makeatletter 
\let\NAT@parse\undefined
\makeatother
\usepackage[colorlinks = true,
            linkcolor = blue,
            urlcolor  = blue,
            citecolor = blue,
            anchorcolor = blue]{hyperref} 

\title{\LARGE \bf
A Game-Theoretic Framework for Network Formation in Large Populations
}

\author{
G\"ok{\c c}e Dayan{\i}kl{\i} and Mathieu Lauri{\`e}re 
\thanks{\textit{Acknowledgments.} Authors thank Agostino Capponi for the discussions and his valuable feedback on the manuscript.\vskip0.1mm Authors' names are alphabetical.}
\thanks{Department of Statistics,
  University of Illinois at Urbana-Champaign, 
  Champaign, IL 61820, USA. 
        {\tt\small gokced@illinos.edu}}%
\thanks{NYU-ECNU Institute of Mathematical Sciences at NYU Shanghai; Shanghai Frontiers Science Center of Artificial Intelligence and Deep Learning; NYU Shanghai, 567 West Yangsi Road, Shanghai, 200126, People’s Republic of China.
        {\tt\small mathieu.lauriere@nyu.edu}}%
}

\begin{document}

\maketitle
\thispagestyle{empty}
\pagestyle{empty}

\begin{abstract}
In this paper, we study a model of network formation in large populations. Each agent can choose the strength of interaction (i.e. connection) with other agents to find a Nash equilibrium. Different from the recently-developed theory of graphon games, here each agent's control depends not only on her own index but also on the index of other agents. After defining the general model of the game, we focus on a special case with piecewise constant graphs and we provide optimality conditions through a system of forward-backward stochastic differential equations. Furthermore, we show the uniqueness and existence results. Finally, we provide numerical experiments to discuss the effects of different model settings.
\end{abstract}

\section{Introduction} 
\label{sec:introduction}

The question of how networks are formed has attracted a growing interest in the past decades, with descriptive models such as the Erd\"os-R\'enyi random graph~\cite{renyi1959random} or the small-world model of~\cite{watts1998collective}, and prescriptive models based on optimization, see e.g.~\cite{bertsekas1991linear}. In many situations, the network is formed as a result of individual interactions between rational agents. Network games model strategic interactions among agents whose payoffs depend on a structured set of relationships, represented as a network. In these games, each agent’s decision, such as investing in relationships, exerting effort, or exchanging information, affects not only their own utility but also that of their neighbors in the network. A central theme is how the network structure influences equilibrium outcomes, including efficiency, stability, and welfare implications. Classical studies on network games include the work of~\cite{jackson2003strategic} on the formation of networks with externalities and~\cite{bala2000noncooperative} on strategic network formation. Applications of network games are diverse, spanning economic markets (e.g., firms forming alliances), social interactions (e.g., opinion dynamics and peer effects), epidemic modeling (e.g., disease spread and vaccination strategies), cybersecurity~\cite{alpcan2010network,omic2009protecting}, and financial systems (e.g., interbank lending and systemic risk)~\cite{capponi2024systemic,acemoglu2015networks}. These models provide insights into the interplay between strategic behavior and network topology, which is crucial for designing policies that promote efficient and resilient networks.

Stochastic differential games extend classical game theory to dynamic settings where agents make decisions in a stochastic environment, often modeled by controlled stochastic differential equations. A key concept in these games is the Nash equilibrium, which characterizes stable strategies where no agent benefits from unilateral deviations. These games are closely related to stochastic optimal control, as each agent's strategy can be viewed as solving an optimal control problem while accounting for the strategic behavior of others. Foundational work in this area includes~\cite{isaacs1969differential} on differential games and~\cite{fleming1989existence} on zero-sum stochastic games and viscosity solutions. Applications of stochastic differential games are widespread, appearing in finance (e.g., portfolio optimization with competition), economics (e.g., dynamic contract theory) and engineering (e.g., control of multi-agent robotic systems). These models provide a rich mathematical framework for analyzing decision-making in uncertain, multi-agent systems. Recently, mean field games~\cite{lasry2007mean,huang2006large} were introduced to study the limit of large population games. Graphon games~\cite{parise2019graphon,caines2021graphon} were introduced to study the limit of many-agent games interacting through a graph structure and used to model applications in epidemics modeling~\cite{aurell2022finite,fabian2022mean} and rumor propagation~\cite{huaning_acc2025} on large networks. However, existing works on graphon games are limited to situations in which the graph structure cannot be changed by the agents, which is not realistic in many applications.

In the present work, we study a model of network formation in the context of stochastic differential games with infinitely many agents. Each agent can choose the strength which she wants to interact with other agents. This leads to a new type of game, in which a key feature is the fact that the control of each agent depends not only on her own index but also on the index of each other agent. This is in contrast with, for example, mean field games and graphon games. We study equilibrium conditions in the form of forward-backward stochastic differential equations. Then, we focus on a linear-quadratic model for which we express the solution in terms of a system of ordinary differential equations. Finally, we show numerical results illustrating the behavior of the solution.

The rest of the paper is organized as follows. In Section~\ref{sec:model}, we present the general network formation model with a continuum of agents. In Section~\ref{sec:speacialmodel}, we consider a special case of piece-wise constant graph. The main results are presented in Section~\ref{sec:main_results}. Numerical algorithm and experiments are discussed in Section~\ref{sec:numerics}. Finally, Section~\ref{sec:conclusion} provides a conclusion and a discussion on future work.

\section{General Model} 
\label{sec:model}

We aim to model the dynamic network formation process among a continuum of agents in a game. We consider a finite time horizon $T > 0$. We will use bold letters to denote functions of time. To simplify the notations,
we will restrict the presentation to one-dimensional states and actions; however, the ideas could be generalized to the multi-dimensional case in a straightforward way. We assume that the number $N$ of noncooperative agents goes to infinity. In the limiting model, each agent is denoted by an index $i\in[0,1]=:I$. At time $t$, agent $i$ has a state $X_t^i \in \RR$ and chooses an interaction (connection) level (or strength) $w_t^i(j)\in[0,1]$ with each agent $j \in I$. In this work, we focus on decisions that are functions of time only. Extensions to state-dependent controls are left for future work. 

For a given control $\bw^i= \big(w^i_t(j)\big)_{t\in[0,T], j\in I}$, agent $i$ is influenced by the aggregate quantity $Z_t^{i,w^i} \in \RR$ defined as:
\small
    \begin{equation*}
        Z_t^{i,\bw^i} = \int_I w_t^i(j)\bar{X}_t^j d\lambda(j),
        \qquad \bar{X}_t^j = \EE[X_t^j],\quad j \in I,
    \end{equation*}
    \normalsize
    where $\lambda$ denotes a measure which is absolutely continuous with respect to the Lebesgue measure on $I$. 
By choosing $w^i_t(j)$, agent $i$ chooses how much she accepts to be influenced by agent $j$'s state. The fact that the aggregate $Z_t^{i,w^i}$ involves $\bar{X}_t^j$ and not $X_t^j$ can be justified in two different ways: either using the exact law of large numbers \cite{sun2006exact,aurell2022stochastic,aurell2022finite}, or viewing each agent as the limit of a sub-population \cite{gao2019graphon,caines2021graphon}.

Given a control profile $(\bw^j)_{j\in I}=:\tilde \bw: I \times I \times [0,T] \to \RR$ used by the population, the aim of agent $i$ is to minimize the following cost over her own control $\bw^i: I \times [0,T] \to \RR$:
\small
\begin{equation*}
    J(\bw^i; \tilde \bw) = \EE\left[ \int_0^T \Big((Z_t^{i, \bw^i}-X_t^i)^2 + \nu^i \int_I (w_t^i(j))^2 d\lambda(j)\Big) dt\right],
\end{equation*}
\normalsize
    where $\nu^i>0$ is an agent-specific constant, subject to the following state dynamics:
    \small
    \begin{equation*}
        \begin{cases}
            dX_t^i &= a^i(Z_t^{i,\bw^i}-X_t^i)dt + \sigma^i dW_t^i
            \\
            d X_t^j &= a^j(Z_t^{j,\tilde \bw^j}-X_t^j)dt + \sigma^j dW_t^j, \qquad j \neq i,
        \end{cases}
    \end{equation*}
    \normalsize
    where $a^i$ and $\sigma^i>0$ are  agent-specific constants.

This model is inspired by the systemic risk model for banks introduced as a mean field game problem in~\cite{carmona2015systemic}. As in that work, the state represents the cash reserve of each agent (i.e., bank)~$i$. The first term in the objective functional reflects the agents’ desire to form connections that minimize their deviation from the aggregate state that is defined as the weighted cash reserve of the connected banks. The second term captures the cost of establishing and maintaining these connections. 
While the model is primarily motivated by the formation of interbank networks, it can also be adapted to other contexts, such as social media network formation, by using alternative objective functions or state dynamics.

The key difference between the current model and the regular graphon game models is that in graphon games, the underlying graph is given exogenously, whereas in the current model, the agents \emph{build} the underlying graphon by choosing the connection levels $\bw$ to find an \textit{equilibrium}.

For example, if agent $i$ observes that her state, $X_t^i$, is too far from the aggregate $Z_t^{i,\bw^i}$ for a certain $\bw^i$, then she could decide to replace $\bw^i$ by $\tilde\bw^i$ chosen such that $Z_t^{i,\tilde\bw^i}$ is closer to $X_t^i$. In other words, agent $i$ can decide to change the her interaction function $\bw^i$ in order to make the aggregate closer to her preference.

In this context, we will look for a Nash equilibrium. We provide the following informal definition before focusing on a special case.
\begin{definition}
    $\hat \bw = (\hat w_t^i)_{t \in [0,T], i \in I}$ is a Nash equilibrium if for every $i \in I$, $\hat \bw^i$ is a minimizer of $J(\cdot; \hat \bw)$. 
\end{definition}
Intuitively, this means that no agent has an incentive to unilaterally change her choice of connection strength $\bw^i$ with other agents.

In the sequel, we focus on a special case for the sake of tractability of the solutions.

\section{Special case: Piecewise-constant graph}
\label{sec:speacialmodel}

In this section, we look at a special case where we model a continuum of agents in $K$ many groups in which the agents in the same group are homogeneous but the agents in different groups are heterogeneous instead of studying a model where there is a continuum of heterogeneous agents. We assume that a proportion $m^k$ of the agents are in group $k$, where $\sum_{k \in \llbracket K \rrbracket} m^k=1$. In this case, we can focus on $K$ many representative agents. Given a control profile $\tilde \bw: \llbracket K \rrbracket \times \llbracket K \rrbracket \times [0,T] \to \RR$, the representative agent in group $k \in \llbracket K \rrbracket:=\{1, \dots, K\}$ aims to minimize the following cost by choosing the interaction level with every other representative agent in each group $\ell$, $\big(w_t^k(\ell)\big)_{t\in[0,T], \ell \in \llbracket K \rrbracket}$:
\small
\begin{equation*}
    J(\bw^k; \tilde \bw)=\EE\left[ \int_0^T \Big((Z_t^k-X_t^k)^2 + \nu^k \sum_{\ell \in\llbracket K \rrbracket} (w_t^k(\ell))^2 m^\ell \Big) dt\right],
\end{equation*}
\normalsize
    where $\nu^k>0$ is an agent-specific constant coefficient and the state dynamics are given as
    \small
    \begin{equation*}
        dX_t^k = a^k(Z_t^k-X_t^k)dt + \sigma^k dW_t^k.
    \end{equation*}
    \normalsize
    Here $a^k$ and volatility $\sigma^k>0$ are exogenous constants and $(W^k_t)_{t\in[0,T]}$ is a Brownian motion that represents the idiosyncratic noise. Furthermore, $Z_t^k$ denotes the aggregate for representative agent $k$ and given as follows:
    \small
    \begin{equation*}
        Z_t^k = \sum_{\ell \in\llbracket K \rrbracket} w_t^k(\ell)\bar{X}_t^\ell m^\ell.
    \end{equation*}
    \normalsize
    where $\bar{X}_t^{\ell} = \EE[X_t^{\ell}]$ and
    \small
    \begin{equation*}
        dX_t^\ell = a^\ell(Z_t^\ell-X_t^\ell)dt + \sigma^\ell dW_t^\ell.
    \end{equation*}  
    \normalsize
     We would like to emphasize that in this special case, the control of representative agent in group $k$ at time $t$, $w_t^k$, is finite ($K$) dimensional instead of infinite dimensional as in Section~\ref{sec:model}. 
\begin{definition}
\label{def:eq_piecewiseconstant}
    $\hat \bw = (\hat w_t^k)_{t \in [0,T], k \in \llbracket K\rrbracket}$ is a Nash equilibrium network formation if for every $k \in \llbracket K\rrbracket$, $\hat \bw^k$ is a minimizer of $J(\cdot; \hat \bw)$. 
\end{definition}

\section{Main Results}
\label{sec:main_results}  

For the simplicity in presentation and notations, we will give the results for the case where there are two groups in the population, i.e., $K=2$; however, the results extend to general $K$ in a straightforward way.
\begin{theorem}
    \label{the:eq_char} \textit{(FBSDE characterization of equilibrium network formation)} Control profile $\hat\bw$ is a Nash equilibrium network formation (see Definition~\ref{def:eq_piecewiseconstant}) if 
\small
\begin{equation}
\label{eq:bestresp_connection}
\hat{w}_t^k(\ell)=\dfrac{\Big[\big(2X_t^k-a^k Y_t^k\big) 2 \bar{X}_t^\ell m^\ell \nu^k m^{-\ell} \Big]}{2\nu^k m^{-\ell}m^{\ell} \Big[\nu^k + (\bar{X}_t^\ell)^2m^{\ell} + (\bar{X}_t^{-\ell})^2 m^{-\ell}\Big]}        
\end{equation}
\normalsize
where $t\in[0,T]$, $k, \ell \in\{1,2\}$, $-\ell$ notation is defined as
\small
\begin{equation*}
    -\ell = \begin{cases}
        1  & \text{if } \ell=2,\\
        2  & \text{if } \ell=1,
    \end{cases}
\end{equation*}   
\normalsize
and where $(\bX, \bY , \tilde \bZ) = \big(X_t^k,Y_t^k,\tilde{Z}_t^k\big)_{k\in\{1,2\}, t\in[0,T]} $ solve the following
forward-backward stochastic differential equation (FBSDE)
system:
\small
\begin{equation*}
 \begin{aligned}
     dX_t^1 &= a^1 (Z_t^1 -X_t^1)dt + \sigma^1 dW_t^1,\\[1mm]
     dX_t^2 &= a^2 (Z_t^2-X_t^2)dt + \sigma^2 dW_t^2,\\[1mm]
     dY_t^1 &= \big(a^1 Y_t^1 + 2(Z_t^1 -X_t^1)\big)dt + \tilde{Z}_t^1dW_t^1,\\[1mm]
     dY_t^2 &= \big(a^2 Y_t^2 + 2(Z_t^2 -X_t^2)\big)dt + \tilde{Z}_t^2 dW_t^2,\\[1mm]
     Z_t^1 &= \hat{w}_t^1(1)\bar{X}_t^1 m^1 + \hat{w}_t^1(2)\bar{X}_t^2 m^2,\\[1mm]
     Z_t^2 &= \hat{w}_t^2(1)\bar{X}_t^1 m^1 + \hat{w}_t^2(2)\bar{X}_t^2 m^2,
 \end{aligned}
 \end{equation*}
 \normalsize
 where $X_0^1\sim \mu_0^1,\ X_0^2\sim \mu_0^2,\ Y_T^1=Y_T^2=0$ and $\bar{X}_t^k = \EE[X_t^k]$ for $k\in\{1,2\}$.
\end{theorem}

The proof of Theorem~\ref{the:eq_char} can be found in the appendix. It relies on an application of stochastic Pontryagin's maximum principle.

\begin{theorem}
\label{thm:existence}
Suppose $T, \nu^k, m^k$ for $k\in\{1,2\}$ are small enough. Then there \textit{exists a unique} solution to the FBSDE system given in Theorem~\ref{the:eq_char}.
\end{theorem}
The proof of Theorem~\ref{thm:existence} can be found in the appendix. It relies on applying Banach fixed point theorem.

\section{Numerical Experiments}
\label{sec:numerics}

\subsection{Algorithm}
In fact, computing the Nash equilibrium does not require solving the above FBSDE: In order to find the expected connection strengths, expected states, and expected aggregates under the Nash equilibrium network formation, it is sufficient to solve a forward backward ordinary differential equation system (FBODE) that characterize the expected values at the Nash equilibrium. To derive this  FBODE, we take the expectation of the FBSDE system in Theorem~\ref{the:eq_char} and we end up with the following equation system.
\small
\begin{equation}
\label{eq:fbode}
 \begin{aligned}
     d\bar{X}_t^1 &= a^1 (\bar{Z}_t^1 -\bar{X}_t^1)dt,\\[1mm]
     d\bar{X}_t^2 &= a^2 (\bar{Z}_t^2-\bar{X}_t^2)dt,\\[1mm]
     d\bar{Y}_t^1 &= \big(a^1 \bar{Y}_t^1 + 2(\bar{Z}_t^1 -\bar{X}_t^1)\big)dt,\\[1mm]
     d\bar{Y}_t^2 &= \big(a^2 \bar{Y}_t^2 + 2(\bar{Z}_t^2 -\bar{X}_t^2)\big)dt,\\[1mm]
     \bar{Z}_t^1 &= \bar{\hat{w}}_t^1(1)\bar{X}_t^1 m^1 + \bar{\hat{w}}_t^1(2)\bar{X}_t^2 m^2,\\[1mm]
     \bar{Z}_t^2 &= \bar{\hat{w}}_t^2(1)\bar{X}_t^1 m^1 + \bar{\hat{w}}_t^2(2)\bar{X}_t^2 m^2,
 \end{aligned}
 \end{equation}
 \normalsize
 where $\bar{X}_0^k= \bar{\mu}_0^k,\ \bar{Y}_T^k=0$ for $k\in\{1,2\}$ and where
 \small
 \begin{equation}
\label{eq:bestresp_connection_expected}
\bar{\hat{w}}_t^k(\ell)=\dfrac{\Big[\big(2\bar{X}_t^k-a^k \bar{Y}_t^k\big) 2 \bar{X}_t^\ell m^\ell \nu^k m^{-\ell} \Big]}{2\nu^k m^{-\ell}m^{\ell} \Big[\nu^k + (\bar{X}_t^\ell)^2m^{\ell} + (\bar{X}_t^{-\ell})^2 m^{-\ell}\Big]}        
\end{equation}
\normalsize
 Above, we used the following notation $\bar{\xi}:=\EE[\xi]$ for the expectation of any random variable $\xi$. 

 We solve this FBODE system by using fixed point algorithm. Convergence of the algorithm follows from standard contraction mapping arguments. The pseudo-code can be found in Algorithm~\ref{algo:Graph}.

{
\begin{algorithm}
\caption{\small Network formation equilibrium algorithm\label{algo:Graph}}

\small
\textbf{Input:} Model parameters: $a^k,\nu^k,m^k$ for $k \in \{1,2\}$; Expected states at $t=0$: $\bar{\mu}_0^k$ for $k \in \{1,2\}$; Time horizon: $T$; Time increments: $\Delta t$, Convergence parameter: $\epsilon$;

\textbf{Output:} Expected equilibrium connection strengths: $\mathbb{E}[\hat{\bw}^k(\ell)] = \bar{\hat{\bw}}^k(\ell)$ for $k, \ell\in \{1,2\}$; Expected States at equilibrium: $\bar{\bX}^k$ for $k \in \{1,2\}$; Expected Aggregates at equilibrium: $\bar{\bZ}^k$ for $k \in \{1,2\}$.

\vskip1mm

\begin{algorithmic}[1]
\STATE Initialize $\bar{\bX}^{k,(0)}=(\bar{X}_0^{k,(0)},\bar{X}_{\Delta t}^{k,(0)},\dots, \bar{X}_{T}^{k,(0)})$ and $\bar{\bY}^{k,(0)}=(\bar{Y}_0^{k,(0)},\bar{Y}_{\Delta t}^{k,(0)},\dots, \bar{Y}_{T}^{k,(0)})$ for $k\in\{1,2\}$.\vskip1mm
\WHILE{$\lVert \bar{\bX}^{1,(j)}- \bar{\bX}^{1,(j-1)}\rVert>\epsilon$ \OR $\lVert \bar{\bX}^{2,(j)}- \bar{\bX}^{2,(j-1)}\rVert>\epsilon$ \OR $\lVert \bar{\bY}^{1,(j)}- \bar{\bY}^{1,(j-1)}\rVert>\epsilon$ \OR $\lVert \bar{\bY}^{2,(j)}- \bar{\bY}^{2,(j-1)}\rVert>\epsilon$ }\vskip1mm
    \STATE{Update $\bar{\bX}^{1, (j+1)}$ by solving the first (forward) ODE in~\eqref{eq:fbode} by using $\bar{\bX}^{2, (j)}$, $\bar{\bY}^{1, (j)}$, and $\bar{\bY}^{2, (j)}$}\vskip1mm
    \STATE{Update $\bar{\bX}^{2, (j+1)}$ by solving the second (forward) ODE in~\eqref{eq:fbode} by using $\bar{\bX}^{1, (j)}$, $\bar{\bY}^{1, (j)}$, and $\bar{\bY}^{2, (j)}$}\vskip1mm
    \STATE{Update $\bar{\bY}^{1, (j+1)}$ by solving the third (backward) ODE in~\eqref{eq:fbode} by using $\bar{\bX}^{1, (j)}$, $\bar{\bX}^{2, (j)}$, and $\bar{\bY}^{2, (j)}$}\vskip1mm
    \STATE{Update $\bar{\bY}^{2, (j+1)}$ by solving the fourth (backward) ODE in~\eqref{eq:fbode} by using $\bar{\bX}^{1, (j)}$, $\bar{\bX}^{2, (j)}$, and $\bar{\bY}^{1, (j)}$}\vskip1mm
\ENDWHILE
\vskip1mm
\STATE{Calculate $\bar{\bw}^k(\ell)$ by plugging in $\bar{\bX}^{1,(j+1)}, \bar{\bX}^{2,(j+1)}, \bar{\bY}^{1,(j+1)}, \bar{\bY}^{2,(j+1)}$ in equation~\eqref{eq:bestresp_connection_expected} for all $k, \ell\in\{1,2\}$.}
\STATE{Calculate $\bar{\bZ}^k$ by plugging in $\bar{\bX}^{1,(j+1)}, \bar{\bX}^{2,(j+1)}, \bar{\bY}^{1,(j+1)}, \bar{\bY}^{2,(j+1)}$ in the last two equations of the system~\eqref{eq:fbode} for all $k, \ell\in\{1,2\}$.}
\RETURN $\bar{\hat{\bw}}^k(\ell)$, $\bar{\bX}^{k,(j+1)}$, $\bar{\bZ}^k$, for $k, \ell \in \{1,2\}$.

\end{algorithmic}
\end{algorithm}}

\subsection{Numerical Results}

To illustrate our results and to understand the effects of model parameters, we provide numerical experiments. We will continue focusing on the setup with two groups. In the plots, group 1 results are shown with blue lines and group 2 results are show red lines.

\begin{figure}[h]
\vskip-2mm
    \centering
    \includegraphics[width=1.0\linewidth]{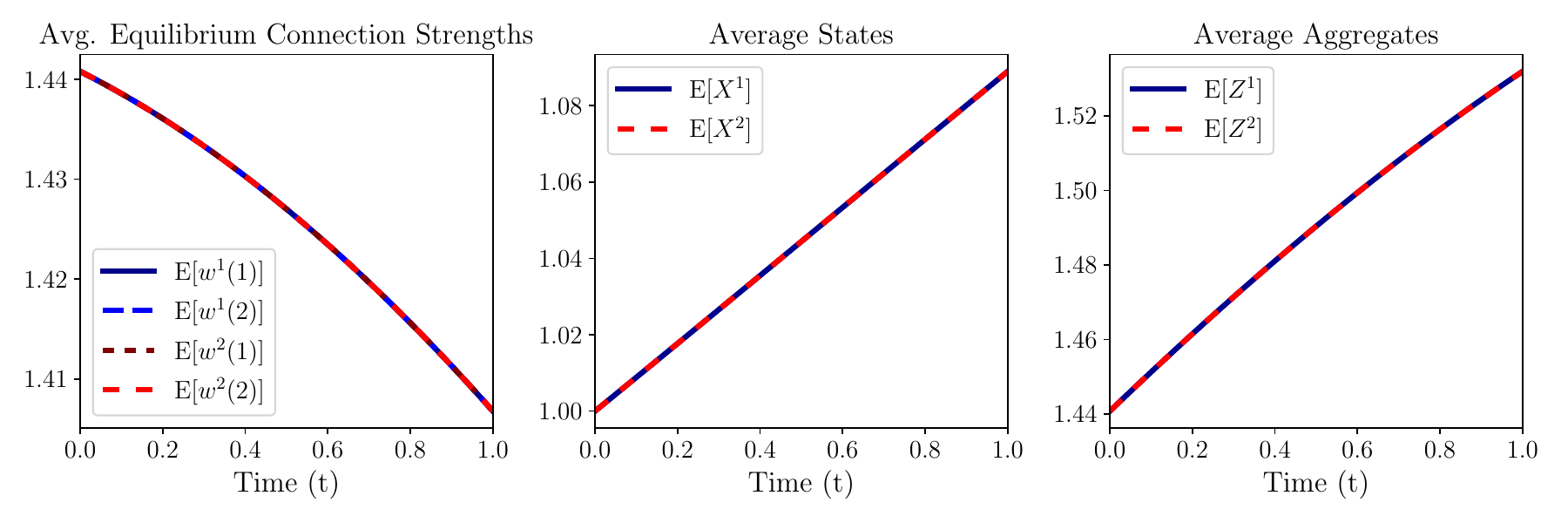}
    \caption{\small \textbf{Base Experiment:} Expected equilibrium connection strengths: $\EE[\hat{w}^k(\ell)]$ for $k, \ell \in\{1,2\}$ (left), Expected States at equilibrium: $\EE[X^k]$ (middle), Expected Aggregates at Equilibrium  $\EE[Z^k]$ for $k\in\{1,2\}$ (right) in the base experiment where the subgroups parameters are equal to each other.}
    \label{fig:EXP0_fullsymm}
    \vskip-3mm
\end{figure}

We will first start with the situation where the proportion of groups are equal to each other ($m^1=m^2=0.5$), and the model parameters are the same for both groups ($a^1=a^2=0.2$, $\nu^1=\nu^2=0.5$, $\bar{\mu}^1_0 =\bar{\mu}^2_0=1.0$). We choose $T=1$ and $\Delta t=0.01$. The parameters set here will constitute our \textit{base experiment parameters}. Since the group parameters are fully symmetric, the agents become fully homogeneous. Therefore, this model will resemble mean field game models and the results of both groups will be equivalent to each other as in Figure~\ref{fig:EXP0_fullsymm}. We want to insist that the results stay the same if proportions of the groups change, since the populations still stay homogeneous due to same parameters.

\begin{figure}[h]
    \centering
    \includegraphics[width=1.0\linewidth]{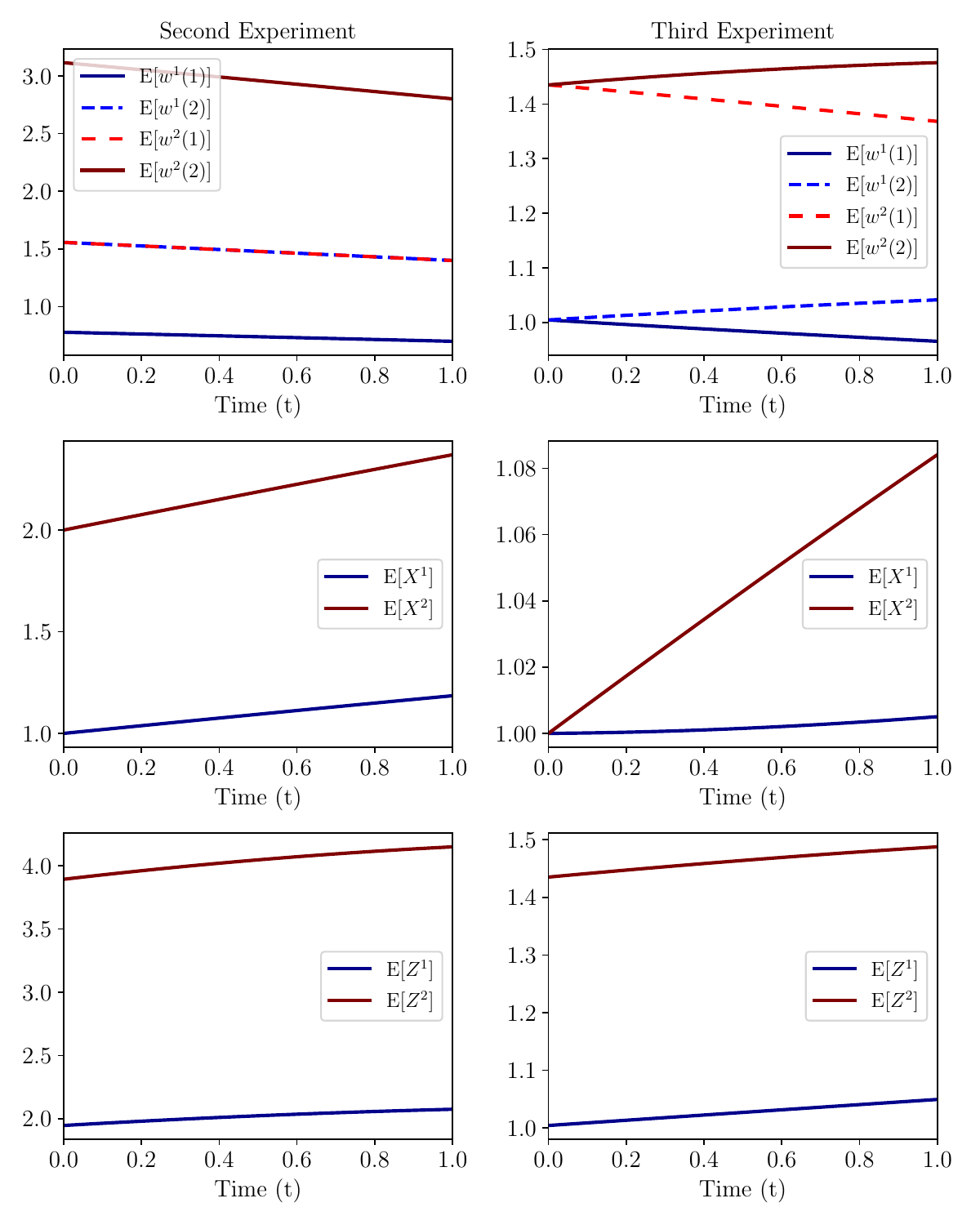}
    \caption{\small \textbf{Second and Third Experiments:} Expected equilibrium connection strengths: $\EE[\hat{w}^k(\ell)]$ for $k, \ell \in\{1,2\}$ (top), Expected States at equilibrium: $\EE[X^k]$ (middle), Expected Aggregates at Equilibrium  $\EE[Z^k]$ for $k\in\{1,2\}$ (bottom) in the second experiment where the initial average states are different: $\bar{\mu}^1_0=1.0$, $\bar{\mu}^2_0=2.0$ (left) and in the third experiment where the cost parameters for different groups are different: ${\nu}^1=1.0$, ${\nu}^2=0.5$ (right).}
    \vskip-5mm
    \label{fig:EXP_diffmu0_diffnu}
\end{figure}

In the \textit{second experiment}, we maintain all conditions from the base experiment, except for the initial expected state values at time 0, which are set to $\bar{\mu}^1_0 = 1.0$ and $\bar{\mu}^2_0 = 2.0$. The results can be seen on the left side of Fig.~\ref{fig:EXP_diffmu0_diffnu}. This change in initial expected state values leads to distinct behaviors among representative agents in each group. Specifically, in the top-left plot of Fig.~\ref{fig:EXP_diffmu0_diffnu}, we observe that the agents in group 2 (with higher initial average state) will have on average higher within group connection strength than the agents in group 1. This behavior helps align the aggregate state more closely with their own. Additionally, we observe that the across group connection strengths of both groups equal to each other due to the fact that other model parameters are the same. In the \textit{third experiment}, we keep all conditions the same as in the base experiment, except for the cost parameters associated with connection strengths, which are set to ${\nu}^1 = 1.0$ and ${\nu}^2 = 0.5$. The results are shown on the right side of Fig.~\ref{fig:EXP_diffmu0_diffnu}. Specifically, in the top-right plot of Fig.~\ref{fig:EXP_diffmu0_diffnu}, we observe that the agents in group 1 (with higher cost parameter for the connection strengths) tend to exhibit on average lower connection strengths (both within the group and across the groups) than the agents in group 2 to offset the higher costs for the same levels of connections. We also see that on average agents in group 1 are less connected within the group than across the group; on the other hand, the agents in group 2 are on average more connected within the group than across the group. Interestingly, we observe that average connections within group 1 gets weaker with time, suggesting that high connection costs may cause group dissolution (i.e., connections getting weaker within the group).

\begin{figure}[h]
    \centering
    \includegraphics[width=1.0\linewidth]{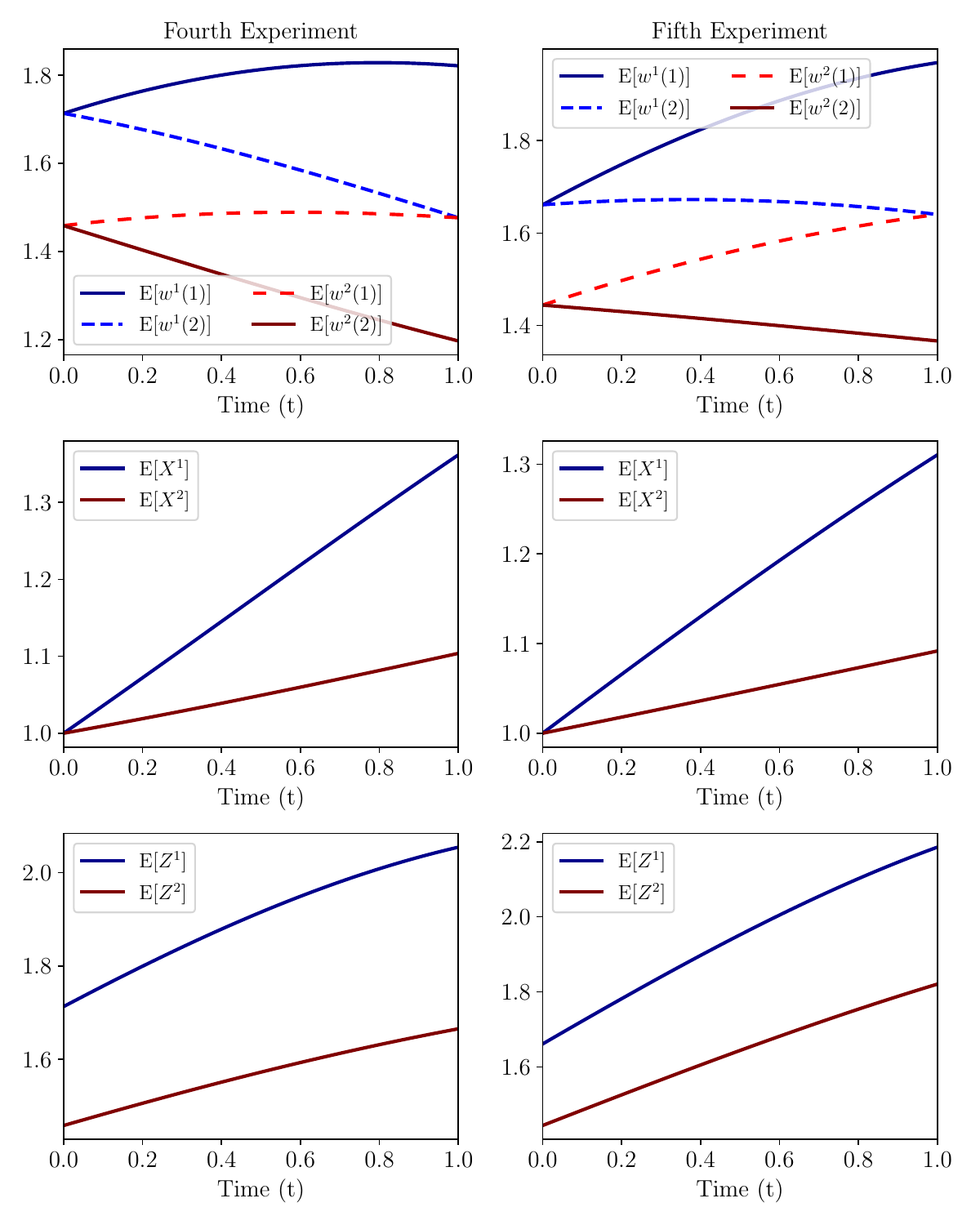}
    \caption{\small \textbf{Fourth and Fifth Experiment:} Expected equilibrium connection strengths: $\EE[\hat{w}^k(\ell)]$ for $k, \ell \in\{1,2\}$ (top), Expected States at equilibrium: $\EE[X^k]$ (middle), Expected Aggregates at Equilibrium  $\EE[Z^k]$ for $k\in\{1,2\}$ (bottom) in the experiment where the drift parameters for different groups are different: ${a}^1=0.5$, ${a}^2=0.2$ (left) and in the experiment where the drift parameters for different groups and the proportions for groups are different: ${a}^1=0.5$, ${a}^2=0.2$ and $m^1=0.1$, $m^2=0.9$ (right).}
    \label{fig:EXP_diffa_diffaprop}
    \vskip-5mm
\end{figure}

In the fourth experiment, we retain all conditions from the base experiment but introduce different drift parameters, setting $a^1 = 0.5$ and $a^2 = 0.2$. This implies that group 1 state process can evolve faster. The results can be seen on the left side of Fig.~\ref{fig:EXP_diffa_diffaprop}. As expected we see that the average state in group 1 increases faster than group 2. This incentivizes group 1 agents to have on average higher levels of aggregates than group 2 agents. Therefore, group 1 agents on average choose higher levels of connection strengths. In the final experiment, we build on the fourth experiment by altering the group proportions, setting $m^1=0.1$, $m^2=0.9$. The results can be seen in Fig.~\ref{fig:EXP_diffa_diffaprop}. Due to lower $m^1$, the aggregates get affected by the average state of group 1 less and this will result in higher connection strengths with group 1 for both groups.   

\section{Conclusion and Future Work}
\label{sec:conclusion}

In this paper, we have proposed a model of network formation for a scenario with infinite population. Each agent decides on the strength which she wants to interact (i.e., connect) with other agents, leading to a differential game where each agent's control is a function of not only her label but also the other agents' labels. For tractable solutions, we focus on a special case. On the theoretical side, we characterize the Nash equilibrium using a forward-backward system of stochastic differential equation reminiscent of McKean-Vlasov FBSDEs. We further give the existence and uniqueness results. We then provide a numerical example to see the effects of different model settings.

Our future work has three main directions. First, we aim to extend our numerical results to incorporate larger number of groups and different model forms for various application domains. Second, we plan to extend the aggregate term to introducing more complex interactions, such as including multiplicative terms as in $Z_t^k 
= \sum_{\ell \in \llbracket L \rrbracket} w_t^k(\ell) w^\ell(k) \bar{X}_t^{\ell} m^{\ell}$. Finally, we intend to analyze more general graph models beyond the piecewise constant.

\bibliographystyle{ieeetr}


\section*{Appendix}

\subsection{Proof of Theorem~\ref{the:eq_char}. }

\begin{proof}
We can write the Hamiltonian of the representative agent in group $k \in \llbracket K \rrbracket$ as follows:
\small
\begin{equation*}
\begin{aligned}
    &H^k(t, w^k, X^k, \bar{X}, Y^k) \\
    &= a^k(\sum_{\ell\in\llbracket K \rrbracket} w^k(\ell)\bar{X}^\ell m^\ell-X^k)Y^k \\
    &+ (\sum_{\ell\in\llbracket K \rrbracket} w^k(\ell)\bar{X}^\ell m^\ell-X^k)^2 + \nu^k \sum_{\ell\in\llbracket K \rrbracket} (w^k(\ell))^2 m^\ell 
\end{aligned}
\end{equation*}
\normalsize
The minimizer of the Hamiltonian is given as 
\small
\begin{equation*}
    \hat{w}_t^k(\ell) = - \frac{a^k m^\ell \bar{X}^\ell_t Y^k_t + 2\bar{X}^\ell_tm^\ell (\sum_{j\in\llbracket K\rrbracket, j\neq \ell} \hat{w}_t^k(j)\bar{X}_t^j m^j - X_t^k)}{2\nu^km^\ell + 2(\bar{X}_t^\ell m^\ell)^2 },
\end{equation*}
\normalsize
for all $k, \ell \in \llbracket K \rrbracket$.

For simplicity in presentation, assume $K=2$. As mentioned before, the results for general $K$ can be obtained in a similar way. Then, we have
\small
\begin{equation}
\begin{aligned}
    \hat{w}_t^1(1) &= -\frac{a^1m^1 Y_t^1\bar{X}_t^1 + 2\bar{X}_t^1m^1 \big(\hat{w}_t^1(2)\bar{X}_t^2m^2 - X_t^1\big)}{2\nu^1 m^1 + 2 (\bar{X}_t^1m^1)^2} \\[1mm]
    \hat{w}_t^1(2) &= -\frac{a^1m^2 Y_t^1\bar{X}_t^2 + 2\bar{X}_t^2m^2 \big(\hat{w}_t^1(1)\bar{X}_t^1m^1 - X_t^1\big)}{2\nu^1 m^2 + 2 (\bar{X}_t^2 m^2)^2} \\[1mm]
    \hat{w}_t^2(1) &= -\frac{a^2 m^1 Y_t^1\bar{X}_t^2 + 2\bar{X}_t^1m^1 \big(\hat{w}_t^2(2)\bar{X}_t^2m^2 - X_t^2\big)}{2\nu^2 m^1 + 2 (\bar{X}_t^1m^1)^2} \\[1mm]
    \hat{w}_t^2(2) &= -\frac{a^2 m^2 Y_t^2\bar{X}_t^2 + 2\bar{X}_t^2m^2 \big(\hat{w}_t^2(1)\bar{X}_t^1m^1 - X_t^2\big)}{2\nu^2 m^2 + 2 (\bar{X}_t^2 m^2)^2} 
\end{aligned}
\end{equation}
\normalsize
We solve the systems for $\hat{w}_t^{k}(\ell)$ for $k, \ell\in \{1, 2\}$:
\small
\begin{equation}
\label{eq:bestresp_connection}
\hat{w}_t^k(\ell)=\dfrac{\Big[\big(2X_t^k-a^k Y_t^k\big) 2 \bar{X}_t^\ell m^\ell \nu^k m^{-\ell} \Big]}{2\nu^k m^{-\ell}m^{\ell} \Big[\nu^k + (\bar{X}_t^\ell)^2m^{\ell} + (\bar{X}_t^{-\ell})^2 m^{-\ell}\Big]}        
\end{equation}
\normalsize
where $-\ell=1$ if $\ell=2$ and $-\ell=2$ if $\ell=1$. We will conclude by using stochastic Pontryagin maximum principle (see e.g.~\cite[Chapter 4]{carmona2018probabilistic}). In this way, the forward equations in the FBSDE system can be written by writing the state dynamics of representative agents in groups $k$ in the equilibrium, i.e., $dX_t^k = a^1 (Z_t^k-X_t^k)dt +\sigma^k dW_t^k$ where $Z_t^k$ is calculated by using the equilibrium connection strengths, i.e., $Z_t^k = \hat{w}^k_t(1) \bar{X}_t^1 m^1 + \hat{w}^k_t(2) \bar{X}_t^2 m^2 $ for all $k \in \{1,2\}$. On the other hand, the backward equations will characterize the derivative of the value functions of representative agents in groups $k$ and will be written as $dY_t^k = -\partial_{x^k} H^k(t, \hat{w}^k, X^k, \bar{X}, Y^k) dt + \tilde{Z}^kdW_t^k$ for $k\in \{1,2\}$ where $\bar{X} = (\bar{X}^1, \bar{X}^1)$ which will result in the following coupled FBSDE system.
 \small
 \begin{equation*}
 \begin{aligned}
     dX_t^1 &= a^1 (Z_t^1 -X_t^1)dt + \sigma^1 dW_t^1,\\
     dX_t^2 &= a^2 (Z_t^2-X_t^2)dt + \sigma^2 dW_t^2,\\
     dY_t^1 &= \big(a^1 Y_t^1 + 2(Z_t^1 -X_t^1)\big)dt + \tilde{Z}_t^1dW_t^1,\\
     dY_t^2 &= \big(a^2 Y_t^2 + 2(Z_t^2 -X_t^2)\big)dt + \tilde{Z}_t^2 dW_t^2,\\
     Z_t^1 &= \hat{w}_t^1(1)\bar{X}_t^1 m^1 + \hat{w}_t^1(2)\bar{X}_t^2 m^2,\\
     Z_t^2 &= \hat{w}_t^2(1)\bar{X}_t^1 m^1 + \hat{w}_t^2(2)\bar{X}_t^2 m^2
 \end{aligned}
 \end{equation*}
 \normalsize
 where $X_0^1\sim \mu_0^1,\ X_0^2\sim \mu_0^2,\ Y_T^1=Y_T^2=0$, and $\hat{w}_t^k(\ell)$ is defined in~\eqref{eq:bestresp_connection}.
\end{proof}

\subsection{Proof of Theorem~\ref{thm:existence}}

\begin{proof}
Taking expectation, we obtain the followings system of ordinary differential equations (ODEs):
\small
\begin{equation*}
 \begin{aligned}
     d\bar{X}_t^1 &= a^1 (\bar{Z}_t^1 - \bar{X}_t^1)dt,
     \\[1mm]
     d\bar{X}_t^2 &= a^2 (\bar{Z}_t^2 - \bar{X}_t^2)dt,
     \\[1mm]
     d\bar{Y}_t^1 &= \big(a^1 \bar{Y}_t^1 + 2(\bar{Z}_t^1 - \bar{X}_t^1)\big)dt,
     \\[1mm]
     d\bar{Y}_t^2 &= \big(a^2 \bar{Y}_t^2 + 2(\bar{Z}_t^2 - \bar{X}_t^2)\big)dt,
     \\[1mm]
     \bar{Z}_t^1 &= \bhatw_t^1(1)\bar{X}_t^1 m^1 + \bhatw_t^1(2)\bar{X}_t^2 m^2,
     \\[1mm]
     \bar{Z}_t^2 &= \bhatw_t^2(1)\bar{X}_t^1 m^1 + \bhatw_t^2(2)\bar{X}_t^2 m^2,
 \end{aligned}
 \end{equation*}
 \normalsize
 where $\bar{X}_0^1 = \EE[X_0^1],$ $\bar{X}_0^2 = \EE[X_0^2],$ $\bar{Y}_T^1=\bar{Y}_T^2=0$ and $\bhatw$ is defined, by taking the expectation in~\eqref{eq:bestresp_connection}, as:
 \small
 \begin{equation*}\bhatw_t^k(\ell)=\dfrac{\Big[\big(2\bar{X}_t^k-a^k \bar{Y}_t^k\big) 2 \bar{X}_t^\ell m^\ell \nu^k m^{-\ell} \Big]}{2\nu^k m^{-\ell}m^{\ell} \Big[\nu^k + (\bar{X}_t^\ell)^2m^{\ell} + (\bar{X}_t^{-\ell})^2 m^{-\ell}\Big]}  .      
\end{equation*}
\normalsize
Let $\varphi^k_x: \RR^4 \to \RR$, $k=1,2$, be defined such that:
\small
\[
    \varphi^k_x(\bar{X}^1_t, \bar{X}^2_t, \bar{Y}^1_t, \bar{Y}^1_t) = a^k (\bar{Z}_t^k - \bar{X}_t^k),
\]
\normalsize
which is the right-hand side in the ODE for $\bar{X}^k_t$. Analogously, let $\varphi^k_y: \RR^4 \to \RR$, $k=1,2$, be defined such that:
\small
\[
    \varphi^k_y(\bar{X}^1_t, \bar{X}^2_t, \bar{Y}^1_t, \bar{Y}^1_t) = \big(a^k \bar{Y}_t^k + 2(\bar{Z}_t^k - \bar{X}_t^k)\big),
\]
\normalsize
which is the right-hand side in the ODE for $\bar{Y}^k_t$.

Notice that $\varphi^k_x$ and $\varphi^k_y$ are locally Lipschitz continuous. Indeed, the only non-linear part comes from the terms with $\bhatw$, which is of the form:
\small
 \begin{equation*}\bhatw_t^k(\ell)=\dfrac{C_1\bar{X}_t^k \bar{X}_t^\ell - C_2 \bar{Y}_t^k \bar{X}_t^\ell }{C_3 + C_4(\bar{X}_t^\ell)^2 + C_5(\bar{X}_t^{-\ell})^2},   
\end{equation*}
\normalsize
where $C_3,C_4,C_5$ are positive constants. 
Let us denote by $L_x^R$ and $L_y^R$ the Lipschitz constants of respectively $\varphi^k_x$ and $\varphi^k_y$ on the ball of radius $R$ (i.e., the set of 4-tuples of continuous functions of time bounded by $R$ in sup-norm).

Now, let $\Phi_x: C([0,T])^2 \to C([0,T])^2$ be the function which maps $(\bar{Y}^1, \bar{Y}^2)$ to the solution $(\bar{X}^1,\bar{X}^2)$ of the first two ODE with given $(\bar{Y}^1, \bar{Y}^2)$.

Let us show that $\Phi_x$ is contractive over a suitable set $\cY$. Let $\cY$ be the subset of $C([0,T])^2$ consisting of $(\bar{Y}^1,\bar{Y}^2)$ which are bounded by constant $C_Y$. Take $(\bar{Y}^1,\bar{Y}^2)$ and $(\bar{Y^\prime}^1,\bar{Y^\prime}^2)$. Denote $\bar{X} = (\bar{X}^1,\bar{X}^2) = \Phi_x(\bar{Y}^1,\bar{Y}^2)$ and $\bar{X^\prime} = (\bar{X^\prime}^1,\bar{X^\prime}^2) = \Phi_x(\bar{Y^\prime}^1,\bar{Y^\prime}^2)$. Let $\tilde{X} = (\tilde{X}^1, \tilde{X}^2) = (\bar{X}^1-\bar{X^\prime}^1,\bar{X}^2-\bar{X^\prime}^2)$ and likewise for $\tilde{Y}$. We have:
\small
\[
    d \tilde{X}^k_t = \left[\varphi^k_x(\bar{X}_t, \bar{Y}_t) - \varphi^k_x(\bar{X^\prime}_t, \bar{Y^\prime}_t)\right] dt
\]
\normalsize
and initial condition $\tilde{X}^k_0 = 0$. Then:
\small
\begin{align*}
    |\tilde{X}^k_t|
    &\le \int_0^t |\varphi^k_x(\bar{X}_s, \bar{Y}_s) - \varphi^k_x(\bar{X^\prime}_s, \bar{Y^\prime}_s)| ds
    \\
    &\le L_x^{C_Y}\int_0^t \left[|\tilde{X}^1_s| + |\tilde{X}^2_s| + |\tilde{Y}^1_s| + |\tilde{Y}^2_s|\right] ds.
\end{align*}
\normalsize
Let $\xi_t = |\tilde{X}^1_t| + |\tilde{X}^2_t|$. Then:
\small
\begin{align*}
    \xi_t
    &\le 2L_x^{C_Y}\int_0^t \xi_s ds + t \beta,
\end{align*}
\normalsize
where $\beta = \max\{\|\tilde{Y}^1\|_\infty, \|\tilde{Y}^2\|_\infty\}$. By Gronwall's lemma,
\small
\begin{align*}
    \xi_t
    &\le 2 t \beta \exp\left(2 L_x^{C_Y} t\right).
\end{align*}
\normalsize
As a consequence, we proved that for all $(\bar{Y}^1,\bar{Y}^2), (\bar{Y^\prime}^1,\bar{Y^\prime}^2) \in \cY$, we have 
\small
\[
\begin{split}
    &\|\Phi_x(\bar{Y}^1,\bar{Y}^2) - \Phi_x(\bar{Y^\prime}^1,\bar{Y^\prime}^2)\|_\infty
    \\
    &\le 2 T  \|(\bar{Y}^1,\bar{Y}^2) - (\bar{Y^\prime}^1,\bar{Y^\prime}^2)\|_\infty \exp\left(2 L_x^{C_Y} T\right).
\end{split}
\]
\normalsize
For $T$ small enough, $\Phi_x$ is a strict contraction. 

We can proceed similarly for the map $\Phi_y: C([0,T])^2 \to C([0,T])^2$, defined as the function which maps $(\bar{X}^1,\bar{X}^2)$ to the solution $(\bar{Y}^1, \bar{Y}^2)$ of the last two ODEs with given $(\bar{X}^1,\bar{X}^2)$. Note that these ODEs are linear so there is no issue of existence. We can show that, for $T$ small enough, $\Phi_y$ is contractive over a suitable set $\cX$. Let $\cX$ be the subset of $C([0,T])^2$ consisting of $(\bar{X}^1,\bar{X}^2)$ which are bounded by constant $C_X$. 

Now, let $\Phi := \Phi_y \circ \Phi_x$. If $T$ small enough and the initial distributions $\mu^1_0,\mu^2_0$ are concentrated close enough to $0$, then $\Phi(\cX) \subset \cX$. Furthermore, by the above argument, $\Phi$ is contractive on $\cX$. By Banach fixed point theorem, $\Phi$ has a unique fixed point in $\cX$. Together with the associated solution to the backward ODEs, they form a solution to the forward-backward ODE system. 
\end{proof}

\end{document}